\documentclass{amsart}
\usepackage{bm}
\usepackage{graphicx, psfrag}

\def\HH{{\mathcal H}}

\def\LL{{\mathcal L}}

\def\bbR{\mathbb{R}}
\def\bbC{\mathbb{C}}
\def\bbZ{\mathbb{Z}}

\def\0{\mathbf{0}}

\def\qedrule{\hbox{\vrule height5.5pt depth0.5pt width3pt}}
\def\qed{\ifmmode \hskip 12pt plus 4pt minus 9pt \qedrule \else
\unskip \nobreak  \hskip 12pt plus 4pt minus 9pt \qedrule \medbreak \fi}

\def\<{\langle}
\def\>{\rangle}

\def\Ga{\Gamma}

\newtheorem{theorem}{Theorem}[section]
\newtheorem{proposition}[theorem]{Proposition}%[section]
\newtheorem{lemma}[theorem]{Lemma}%[theorem]
\newtheorem{corollary}[theorem]{Corollary}%[theorem]
\newtheorem{conjecture}[theorem]{Conjecture}%[theorem]

\theoremstyle{definition}
\newtheorem{remark}[theorem]{Remark}%[theorem]
\newenvironment{tpmatrix}{\left(\begin{smallmatrix}}{\end{smallmatrix}\right)}

\numberwithin{equation}{section}

\begin{document}
\title{Extensions of positive definite functions on amenable groups}


\author{M. Bakonyi and D. Timotin}
\address{Department of Mathematics and Statistics, Georgia State University, P.O. Box 4110, Atlanta, GA 30302-4110, USA}
\email{mbakonyi@gsu.edu}
\address{Institute of Mathematics of the Romanian Academy, PO Box 1-764, Bucharest 014700, Romania} \email{Dan.Timotin@imar.ro}
%\date{\today}
\thanks{This work has been partially supported by NSF
Grant 12-21-11220-N65 and by Grant 2-Cex06-11-34/2006 of the Romanian Government. A visit of the first author to Bucharest has been supported by SOFTWIN}

\date{}

\begin{abstract}
Let $S$ be a subset of a amenable group $G$ such that $e\in S$ and $S^{-1}=S$. The main result
of the paper states that if the Cayley graph of $G$ with respect to $S$ has a certain combinatorial
property, then every positive definite operator-valued function on $S$ can be extended to a 
positive definite function on $G$. Several known extension results are obtained as a corollary.
New applications are also presented. 
\end{abstract}

\maketitle

\section{Introduction}

Let $G$ be a group.  A function $\Phi :G\to {\mathcal L}({\mathcal H})$ is called 
{\it positive definite} if for every $g_1,...,g_n\in G$ the operator matrix
$\{\Phi (g_i^{-1}g_j\}_{i,j=1}^n$ is positive semidefinite. Let $S\subset G$
be a \emph{symmetric} set; that is, $e\in S$ and $S^{-1}=S$. A function $\phi :S\to {\mathcal L}({\mathcal H})$ is called 
{\it (partially) positive definite} is for every $g_1,...,g_n\in G$ 
such that $g_i^{-1}g_j\in S$ for all $i,j=1,...,n$, 
$\{\phi (g_i^{-1}g_j\}_{i,j=1}^n$ is a positive semidefinite operator matrix.
Extensions of positive definite functions on groups have a long history, starting
with the Trigonometric Moment Problem of Carath\' eodory and Fej\' er and 
Krein's Extension Theorem. Recently, it has been proved in \cite{B1} that every positive
definite operator-valued function on a symmetric interval in an ordered
abelian group can be extended to a positive definite function on the whole group.  By
different techniques, the same extension property was shown to be true in \cite{BT} for functions
defined on words of length $\le m$ in the free group with $n$ generators. In the present paper
we extend the result to a class of subsets of amenable groups which satisfy a certain 
combinatorial condition. The result turns out to be more general than the main
result in \cite{B1} and it is obtained by much simpler means. Our main result was also
influenced by \cite{Exel}, where a version of Nehari's Problem was solved for
operator functions on totally ordered amenable groups.  

Let $G$ be a locally compact group. A
{\it right invariant mean} 
$m$ on $G$ is a state on $L^\infty(G)$ which satisfies 

\[
m(f)=m(f_x),
\]
for all $x\in G$,
where $f_x(y)=f(yx)$. In case there exists a right invariant mean on $G$, $G$ is
called {\it amenable}. We will occasionally write $m^x(f(x))$ for $m(f)$. There exist many other equivalent characterizations of amenability
\cite{Dix}.

For graph theoretical notions we refer the reader to \cite{Gol}.
By a {\it graph} we mean a pair $G=(V,E)$ in which $V$ is a set called the {\it vertex set} and $E$ is 
a symmetric nonreflexive binary relation on $V$, called the {\it edge set}. We consider in general the vertex set
to be infinite. A graph is called {\it chordal} if every finite simple cycle $[v_1,v_2,...,v_n,v_1]$
in $E$ with $n\ge 4$ contains a chord, i.e. an edge connecting two nonconsecutive
vertices of the cycle.  Chordal graphs play an important role
in the extension theory of positive definite matrices (\cite{GJSW} and \cite{T}).
 
Let $G$ be a group. If $S\subset G$ is symmetric,
we define the {\it Cayley graph} of $G$ with respect to $S$ (denoted $\Gamma(G, S)$)
as the graph whose vertices are the elements of $G$, while $\{x,y\}$ is an edge
iff $x^{-1}y\in S$.

\section{The main result}

The basic result of the paper is the following.

\begin{theorem}\label{th:main}
Suppose $G$ is amenable, and $S\subset G$.
If $\Gamma(G,S)$ is chordal,
then any positive definite function $\phi$ on $S$ admits a positive 
definite extension $\Phi$ on $G$.
%
%(ii) If $\Gamma(G,S)$ is not chordal, then there exists positive definite functions
%on $S$ that do not have positive definite extensions.
\end{theorem}

\begin{proof}  
Consider the partially positive semidefinite kernel $k:G\times G\to \LL(\HH)$, defined 
only for pairs $(x,y)$ for which $x^{-1}y\in S$, by the formula
\[
k(x,y)=\phi(x^{-1}y).
\]
Since the pattern of specified values for this kernel is chordal by assumption,
it follows from~\cite{T} that $k$ can be extended to a positive semidefinite
kernel $K:G\times G\to\LL(\HH)$. Note that $K(x,y)$ has no reason to
depend only on $x^{-1}y$.

For any $x, y\in G$, the operator matrix 
$\begin{tpmatrix}
\phi(e) & K(x,y)\\ K(x,y)^* & \phi(e)
\end{tpmatrix}$ 
is positive semidefinite, whence it follows that $K(x,y)^*K(x,y)\le \phi(e)^2$.
In particular, all operators $K(x,y)$, $x,y\in G$, are bounded by a common
constant.

Fix then $\xi, \eta\in \HH$ and $x\in G$. The function $F_{x; \xi, \eta}:G\to\bbC$,
defined by 
\begin{equation}
	\label{eq:215}
F_{x; \xi, \eta}(y)= \< K(yx, y) \xi, \eta\>
\end{equation}
is in $L^\infty(G)$. Define then $\Phi:G\to \LL(\HH)$ by
\begin{equation}
	\label{eq:216}
\<\Phi(x)\xi, \eta\>=m(F_{x; \xi, \eta}). 
\end{equation}

We claim that $\Phi$ is a positive definite function. Indeed, take arbitrary vectors $\xi_1,\dots, \xi_n\in\HH$. We have
\[
\sum_{i,j=1}^n\< \Phi(g_i^{-1}g_j)\xi_i, \xi_j\>=
\sum_{i,j=1}^n m(F_{g_i^{-1}g_j;\xi_i,\xi_j})
= \sum_{i,j=1}^n m^y\left(\<K(y g_i^{-1}g_j, y)\xi_i, \xi_j\>\right).
\]
Take one of the terms in the last sum; the mean $m$ is applied to the function $y\mapsto \<K(y g_i^{-1}g_j, y)\xi_i, \xi_j\>$. The right invariance of $m$ implies that we may apply the change of variable $z=yg_i^{-1}, \ y=zg_i$, and thus
\[
 m^y\left(\<K(y g_i^{-1}g_j, y)\xi_i, \xi_j\>\right)=
 m^z\left(\<K(z g_j, g_i z)\xi_i, \xi_j\>\right).
\]
Therefore
\[
\sum_{i,j=1}^n\< \Phi(g_i^{-1}g_j)\xi_i, \xi_j\>=
\sum_{i,j=1}^n m\big(\<K(z g_j, g_i z)\xi_i, \xi_j\>\big)=
m\Big( \sum_{i,j=1}^n\<K(z g_j, g_i z)\xi_i, \xi_j\>\Big).
\]
But the positivity of $K$ implies that, for each $z\in G$,
\[
\sum_{i,j=1}^n\<K(z g_j, g_i z)\xi_i, \xi_j\>\ge 0.
\]
Since $m$ is a positive functional, it follows that indeed $\Phi$ is positive definite.
On the other hand, for $x\in S$, the function $F_{x; \xi, \eta}$ is constant, equal to $\<\phi(x)\xi, \eta\>$. Therefore 
$\Phi$ is indeed the desired extension of~$\phi$.
\end{proof}

%\smallskip
%(ii) Suppose $\Gamma(G,S)$ is not chordal; it has then a minimal cycle
%of length at least 4, say $[a_1, \dots , a_n]$. In particular, $x_i=a_i^{-1}a_{i+1}\in S$
%for $1\le i \le n-1$ and $x_n=a_n^{-1}a_1\in S$.
%
%Define then $\phi:S\to\bbC$ by $\phi(e)=1$, $\phi(x_i)=1$ for $1\le i\le n-1$, and $\phi(x)=0$
%for all other $x\in S$.
%
%We claim that $\phi$ is positive definite on~$S$. We have to check that
%for any clique $C$ in $\Gamma(G,S)$ the matrix
%$M_C=(\phi(y^{-1}y')_{y,y'\in C}$ is positive semidefinite. 
%If we consider the pairs $\{y, y'\}$ such that $y^{-1}y'=x_i$ for some $i\ge 1$, then
%the fact that $\gamma$ is a cycle implies that these pairs are all mutually disjoint.
%Therefore the matrix $M_C$ is a direct sum of the identity matrix and a certain number
%of copies of $\begin{tpmatrix} 1 &1\\1 &1 \end{tpmatrix}$, and is therefore positive semidefinite.
%
%
%On the other hand, suppose $\Phi$ is a positive definite extension of $\phi$. Consider
%the matrix $P=(\Phi(a_i^{-1}a_j)_{i,j=1}^n$. Since $\phi$ is an extension of $\Phi$,
%$P$ has  on the main diagonal, as well on the diagonals above and below. Using the fact
%that the matrix 
%\[
%\begin{pmatrix}
%1 &1& z\\ 1&1&1 \\ z&1&1
%\end{pmatrix}
%\]
%is positive definite only for $z=1$, one shows immediately by recurrence that all 
%elements of $P$ must be 1. But $\Phi(a_n^{-1}a_1)= \phi(a_n^{-1}a_1)=\phi(x_n)=0$. The contradiction obtained shows that $\phi$ has no positive
%definite extension.
%\end{proof}

\begin{remark}
\label{rem1}
The chordality of $\Gamma (G,S)$ means that
for every finite cycle $[g_1,...,g_n,g_1]$, $n\ge 4$, at least one $\{g_i, g_{i+2}\}$ 
(with $g_{n+1}=g_1$ and $g_{n+2}=g_2$) is an edge. Denoting $\xi _{k}=g_kg_{k+1}^{-1}$,  
the condition is equivalent to:
$\xi _1,...,,\xi _n\in S$, 
$\xi _1\xi_2 \cdots \xi _n=e$, $n\ge 4$, implies that there 
exist $i=1,...,m$ such that $\xi _i\xi _{i+1}\in S$ (here $\xi _{n+1}=\xi _1$).
\end{remark}

\begin{remark}
\label{rem2}
Let $\Lambda \subset G$ be such that $e\in \Lambda$, and $e$ cannot be written as a
product of elements in $\Lambda $ different from $e$, and let $S=\Lambda \Lambda ^{-1}$.
Assume we have that $S=
\Lambda \cup \Lambda ^{-1}$. Then $\xi _1\xi _2\cdots \xi _n=e$, with $\xi _1,...,\xi _n
\in S$, implies the existence of $k$ such that $\xi _k\in \Lambda$ and $\xi _{k+1}\in
\Lambda ^{-1}$, thus $\xi _k\xi _{k+1}\in S$, implying $\Gamma (G,S)$ is chordal. 
\end{remark}

We conjecture the following reciprocal of Theorem 
\ref{th:main}.

\begin{conjecture}
\label{conj}
For every $S\subset G$ such that $\Gamma (G,S)$
is not chordal there exists a positive definite function $\phi :S\to \LL(\HH)$
which does admit a positive definite extension to $G$. 
\end{conjecture}

The following examples strongly suggest that the above conjecture has a positive answer.
Let $G=\bbZ^2$
and let $S=\bbZ^2-\{(1,1), (-1,-1)\}$, 
the minimal number of points that can be excluded. 
Then $(0,0)$, $(0,1)$, $1,1)$, and
$(-1,0)$ form a chordless cycle of length $4$ in $\Gamma (G,S)$. Define $\phi :S\to
M_2({\mathbb C})$ by $\phi ((0,0))=\begin{tpmatrix}
                                1 & 0\\
				0 & 1\\
                               \end{tpmatrix}
$, $\phi ((1,0))=\begin{tpmatrix}
             0 & 0\\
	     1 & 0\\
	                  \end{tpmatrix}
			  $, $\phi ((0,1))=
			  \begin{tpmatrix}
			  0 & 1\\
			  0 & 0\\
			  \end{tpmatrix}
			  $, and $\phi (g')=0$, for every $g'\in S- \{
			  (0,0), (1,0), (-1,0), (0,1), (0,-1)\}$.
Let $K$ be a maximal clique of $\Gamma (G,S)$. We may
assume that $(0,0)\in K$, in which case $(1,1)\not\in K$. This fact
implies that the matrix $\{\phi (x-y)\}_{x,y\in K}$ can be written as a direct 
sum of copies of $\begin{tpmatrix}
1 & 0\\
0 & 1\\
                  \end{tpmatrix}
$, and $\begin{tpmatrix}
1 & 1\\
1 & 1\\
        \end{tpmatrix}
$, so $\phi $ is positive definite. Assume that $\phi $ admits
a positive definite extension $\Phi $ to $G$. Then, 
since 
\[
\begin{pmatrix}
\Phi ((0,0)) & \Phi ((1,0))^* & \Phi ((1,1))^*\\
\Phi ((1,0)) & \Phi ((0,0)) & \Phi ((0,1))^*\\
\Phi ((1,1)) & \Phi ((0,1)) & \Phi ((0,0))\\
 \end{pmatrix} \ge 0
 \]
and 
  \[
 \begin{pmatrix}
\Phi ((0,0)) & \Phi ((0,1))^* & \Phi ((1,1))^*\\
\Phi ((0,1)) & \Phi ((0,0)) & \Phi ((1,0))^*\\
\Phi ((1,1)) & \Phi ((1,0)) & \Phi ((0,0))\\
 \end{pmatrix} \ge 0,
 \]
  it follows that $\Phi ((1,1))=
 \begin{tpmatrix}
 1 & 0\\
 0 & 1\\
 \end{tpmatrix}
 $. Since 
 \[
\begin{pmatrix}
\Phi ((0,0)) & \Phi ((1,1))^* & \Phi ((2,1)))^*\\
\Phi ((1,1)) & \Phi ((0,0)) & \Phi ((1,1))^*\\
\Phi ((2,1)) & \Phi ((1,1)) & \Phi ((0,0))\\
 \end{pmatrix} \ge 0
 \]
  the $(2,1)$ entry of $\Phi((2,1))$ equals $1$,
  contradicting the fact that  $\Phi((2,1))=\phi ((2,1))=0$.
  This implies that $\phi $ does not admit a positive definite extension
  to $\bbZ^2$. 
  
  \bigskip

  Let $\Lambda \subset {\mathbb Z}^d$ be a finite set. By the definition introduced in \cite{Ru}, a sequence
  $\{c_k\}_{k\in \Lambda -\Lambda }$ of complex numbers is called {\it positive definite with respect to $\Lambda $} if 
  the matrix $\{c_{k-l}\}_{k,l\in \Lambda}$ is positive definite. This definition is weaker then the one
  used in this paper, since it requires only a single matrix built on the given data to be positive definite.
  A finite subset $\Lambda \subset {\mathbb Z}^d$ is said to posses the {\it extension property} if every 
   sequence $\{c_k\}_{k\in \Lambda -\Lambda }$ admits a positive extension to ${\mathbb Z}^d$. Let $R(0,n)=
   \{0\}\times \{0,1,...,n\}$, $R(1,n)=\{0,1\}\times \{0,1,...,n\}$, and $S(1,n)=R(1,n)-\{(1,n)\}$. The following 
   is the main result of \cite{BN2}.
   \begin{theorem}
   \label{bn}
   A finite $\Lambda \subset {\mathbb Z}^2$ has the extension property if and only if $\Lambda $ is the translation
   by a vector in ${\mathbb Z}^2$ of a set isomorphic to one of the following sets: $R(0,n)$, $R(1,n)$, or
   $S(1,n)$, $n\ge 0$.
   \end{theorem}
   
   Let $\Lambda =R(1,n)$, when $S=\Lambda -\Lambda =\{-1,0,1\}\times \{-n,...,0,...,n\}$. By the 
   previous theorem, every scalar positive definite sequence with respect to $\Lambda $ on $S$ admits
   a positive definite extension to ${\mathbb Z}^2$. The points $(0,0)$, $(-1,n)$, $(0, 2n)$, and $(1,n)$ 
   form a chordless cycle in $\Gamma ({\mathbb Z}^2, S)$, and for every Hilbert space ${\mathcal H}$
with ${\rm dim}{\mathcal H}\ge 2$, there exists a sequence   
  $\{C_k\}_{k\in S}$ of operators on ${\mathcal H}$ which is
positive definite (in the stronger sense), but does not admit a positive definite extension to ${\mathbb Z}^2$.
The same is true for the sets $S(1,n)$ as well. We will present next the details 
concerning the different behaviour of scalar and operator sequences for a subset of ${\mathbb Z}^2$ not covered by Theorem \ref{bn}.
                        
Let $G=\bbZ^2$,  
$m,n\in {\mathbb N}$, $m,n\ge 2$,  and let $S$
consist of the points $(k,0)$, $|k|\le m$ together
with the points $(0,l)$, $|l|\le n$. Let $\{C_{kl}\}_{(k,l)\in S}$ be a
positive definite sequence of operators. The positive
definiteness condition is equivalent to 

\begin{equation}
\label{toep}
\begin{pmatrix}
C_{00} & C_{10}^* & \cdots & C_{m0}^*\\
C_{10} & C_{00} & \cdots & C_{m-1,0}^*\\
\vdots & \ddots & \ddots & \vdots\\
C_{m0} & C_{m-1,0} & \cdots &C_{00}\\
 \end{pmatrix} \ge 0 
 \end{equation}
 and
 \begin{equation}
 \label{toep1}
 \begin{pmatrix}
C_{00} & C_{01}^* & \cdots & C_{0n}^*\\
C_{01} & C_{00} & \cdots & C_{0,n-1}^*\\
\vdots & \ddots & \ddots & \vdots\\
C_{0n} & C_{0,n-1} & \cdots & C_{00}\\
  \end{pmatrix} \ge 0.
  \end{equation}
 In case $\{c_{kl}\}_{(k,l)\in S}$ is the sequence defined by 
 $c_{k0}=e^{ik\alpha }$ and $c_{0l}=e^{il\beta }$, the matrices
 in (\ref{toep}) are rank $1$ positive definite Toeplitz matrices and 
 $c_{kl}=e^{ik\alpha }e^{il\beta }$, $(k,l)\in \bbZ^2$ is a positive
 definite extension to $\bbZ^2$ of the initial sequence. It is a classical
 result of Carath\' eodory and Fej\' er that every positive definite Toeplitz matrix
 is a positive linear combination of rank $1$ positive definite Toeplitz
 matrices. This implies that the positive semidefiniteness of the matrices
 in (\ref{toep}) guarantees the existence of a positive definite extension
 to $\bbZ^2$ of every positive definite sequence
 $\{c_{kl}\}_{(k,l)\in S}$ of complex numbers. 
 
 Let $U_1$ and $U_2$ be two 
 noncommuting unitary operators on a Hilbert space ${\mathcal H}$ with
 ${\rm dim}{\mathcal H}\ge 2$. Defining $C_{00}=I$, $C_{k0}=U_1^k$, and $C_{0l}=U_2^l$, 
 the matrices in (\ref{toep}) and (\ref{toep1})
 are positive semidefinite. Assuming the sequence $\{C_{kl}\}_{(k,l)\in S}$ admits a positive 
 definite extension to $\bbZ^2$, the operator $C_{11}$ has to
 simultaneously verify the conditions
 $\begin{pmatrix}
 C_{00} & C_{01}^* & C_{11}^*\\
 C_{01} & C_{00} & C_{10}^*\\
 C_{11} & C_{10} & C_{00}\\
 \end{pmatrix} \ge 0
 $ and $
 \begin{pmatrix}
 C_{00} & C_{10}^* & C_{11}^*\\
 C_{10} & C_{00} & C_{01}^*\\
 C_{11} & C_{01} & C_{00}\\
 \end{pmatrix} \ge 0
 $. For our data, the above conditions are equivalent to $C_{11}=U_2U_1$,
 respectively $C_{11}=U_2U_1$, which is false, since $U_1$ and $U_2$ 
 do not commute. Thus $\{C_{kl}\}_{(k,l)\in S}$ does not
 admit any positive definite extension to $\bbZ^2$. 

\begin{proposition}
\label{lulu}
Let $0\in S=-S$ be a finite subset of $\bbZ^2$ such that $\Gamma (\bbZ^2, S)$ is chordal and 
$S$ spans $\bbZ^2$. Then
$S$ is infinite.
\end{proposition}
 
\begin{proof}
Suppose $S\subset \bbZ^2$ is finite and $\Gamma (\bbZ^2, S)$
is chordal. There are a finite number of directions among the elements of $S$; suppose 
the elements of maximum length in each of these directions, together with their inverses, 
are enumerated $s_1,s_2,\ldots, s_{2n}$ in the order of their arguments. 

For a positive integer $N$ consider  the cycle $[x_0, x_2, \dots x_{2nN-1}, x_0]$ in $\Gamma (\bbZ^2, S)$, defined as follows: $x_0=0$, $x_k-x_{k-1}=s_j$ if $(j-1)N< k \le jN$.
We claim that, if $N$ sufficiently large, this is a cycle with no chords.

Indeed, suppose $\{x_k,x_l\}$ is an edge with $l-k\ge 2$.
The points $x_0,\dots, x_{2nN-1}$ form a polygon $P$ with $2n$ sides  $A_j$ parallel to $s_j$ respectively, each side containing $N$ points $x_k$. We have the following possibilities:

---If $x_k$ and $x_l$ are on the same side $A_j$ of $P$, then $x_l-x_k= (l-k)s_j$ would be an element of $S$ colinear with $s_j$, but longer, which is not possible. 

---If $x_k\in A_j$, $x_l\in A_{j+1}$, then the argument of $x_l-x_k$ would be strictly between the arguments of $s_j$ and $s_{j+1}$: again a contradiction.

---Finally, we may chose $N$ sufficiently large such that, if $x_k$ and $x_l$ are on nonconsecutive sides of $P$, then $x_l-x_k$ has length larger than any element of $S$.

So the cycle obtained has no chords, contrary to the chordality assumption in the hypothesis. Thus $S$ must be infinite.
\end{proof}

\begin{remark}
If Conjecture \ref{conj} is true, then
Lemma \ref{lulu} would imply that for every  
finite $S\subset \bbZ^2$ such that $0\in S=-S$ and $S$ spans $\bbZ^2$,
there exists a positive definite function on $S$ which
does not admit a positive definite extension to $\bbZ^2$.
\end{remark}

\section{Applications}

\subsection{Ordered groups and related questions}

Suppose $G$ is a (left or right) totally ordered group. Take $a\in G$, $a\ge e$, and define
$\Lambda =[e,a)$, and
$S=(a^{-1}, a)$. Then $e$ cannot be written as a product of
elements in $\Lambda $ and $S=\Lambda \Lambda ^{-1}=\Lambda \cup \Lambda ^{-1}$.   
Then by Remark \ref{rem2} the graph $\Gamma(G,S)$ is chordal.
%Indeed, if $[x_1, \dots, x_n]$  is a cycle in $\Gamma(G,S)$, 
%suppose $x_1=\max x_i$, and, for instance, $x_n\ge x_2$. Then 
%\[
%e=x_2^{-1}x_2\le x_2^{-1} x_n \le x_2^{-1} x_1 \le a. 
%\] 
%Therefore $\{x_2, x_n\}$ is an edge in $\Gamma(G,S)$, and therefore 
%there are no minimal cycles of length $>3$. The hypothesis of Theorem~\ref{th:main}
%is therefore satisfied;
Thus, in an amenable totally ordered group any 
positive definite function defined on a symmetric interval can be extended to
the whole group.

The same argument yields the following more general result.

\begin{proposition}\label{pr:ordered}
Suppose $G$ is amenable, while $G'$ is a totally ordered group, with unit $e'$.
Let $g:G\to G'$ be a group morphism. Take $a'\in G'$, $a'\ge e'$, and $S=g^{-1}((a'{}^{-1}, a'))$.
Then any positive definite operator function on $S$ can be extended to a positive definite
function on the whole group.
\end{proposition}

The above proposition has the following consequence which represents the main result of
\cite{B1}. The proof derived here is much simpler.

\begin{corollary}
\label{cucu}
Let $G_1$ be a totally ordered abelian group, $a\in G_1$, $a>0$, and let $G_2$ be an abelian group. Then any
positive definite operator function on $(-a,a)\times G_2$ can be extended to a positive definite
function on $G_1\times G_2$. 
\end{corollary}

Several well-known results, such that the Classical Trigonometric Moment Problem and Krein's Extension
Theorem are particular cases of Corollary \ref{cucu}. 
Another simple application of Corollary  \ref{cucu} is the following. Take
$\alpha, \beta \in \bbR$, and define $g:\bbZ^2\to\bbR$ by $g(m, n)= \alpha m+\beta n$.
Thus, all positive definite functions defined on the strip $|\alpha m+\beta n|< a$
can be extended to a positive definite function on $\bbZ^2$.

A more interesting example for Proposition \ref{pr:ordered} is given by the Heisenberg group $H$ over the integers. This is
the group of matrices  of the form
\[
X_{m,n,p}=\{
\begin{tpmatrix}
1 & m & p\\ 0 & 1 &n\\ 0& 0&1
\end{tpmatrix}, \; m,n,p\in \bbZ \}.
\]
It is an amenable group, and for any $\alpha ,\beta \in \bbR$,
we can consider the morphism $g:H\to \bbR$,
given by $g(X_{m,n,p})=\alpha m+\beta n$. Thus any positive definite function defined 
on the set $\{X_{m,n,p}: |\alpha m+\beta n|< a\}$
 can be extended to a positive definite function on $H$.

\subsection{Trees and Cayley graphs}

For this application we need some supplementary preliminaries.
If $\Gamma=(V,E)$ is a graph,  the distance $d(v,w)$
between two vertices is defined as
\[
d(v,w)=\min\{n: \exists v=v_0, v_1, \dots , v_n=w, \mbox{ such that } \{v_i,v_{i+1}\}
\in E(\Gamma)   \}.
\]
We define the graph $\hat{\Gamma}_n$ that has the same vertices as $\Gamma$,
while $\{v,w\}$ is an edge of $\hat{\Gamma}_n$ if and only if $d(v,w)\le n$.

A graph without any simple cycle is called a {\it tree}. If $x$ and $y$ are two
distinct vertices of a tree, then $P(x,y)$ denotes the unique simple path joining $x$
and $y$.

\begin{lemma}\label{le:chordal}
If $\Gamma$ is a tree, then $\hat\Gamma_n$ is chordal for any $n\ge 1$.
\end{lemma}

\begin{proof} Take a minimal cycle $C$ of length $>3$ in $\hat \Gamma_n$. Suppose $x,y\in C$
maximize the distance between any two points of $C$. If $d(x,y)\le n$, then
$C$ is a clique, which is a contradiction. Thus $x$ and $y$ are
not adjacent in $\hat \Gamma_n$. Suppose $v,w$ are the two vertices of $\hat \Gamma_n$ adjacent
to $x$ in the cycle $C$. Now $P(x,v)$ has to pass through a vertex
which is on $P(x,y)$, since otherwise the union of these two paths
would be the minimal path connecting $y$ and $v$, and it would
have length strictly larger than $d(x,y)$. Denote by $v_0$ the
element of $P(x,v)\cap P(x,y)$ which has the largest distance to
$x$; since $d(y,v)=d(y,v_0)+d(v_0,v)\le d(y,x)=d(y,v_0)+d(v_0,x)$,
it follows that $d(v_0,v)\le d(v_0,x)$.

Similarly, if $w_0$ is the element of $P(x,w)\cap P(x,y)$ which has
the largest distance to $x$, it follows that $d(w_0,w)\le
d(w_0,x)$.

Suppose now that $d(v_0,x)\le d(w_0,x)$. Then
\begin{align*}
d(v,w)&= d(v,v_0)+d(v_0,w_0)+d(w_0,w)\\
&\le d(x, v_0)+d(v_0,w_0)+d(w_0,w)=d(x,w)\le n,
\end{align*}
since $w$ is adjacent to $x$. Then $(v,w)\in E$, and $C$ is not
minimal: a contradiction. Thus $\hat \Gamma_n$ is chordal.
\end{proof}

%%%%

It is worth mentioning that $\Gamma$ chordal does not necessarily
imply $\hat{\Gamma}_n$ chordal. For instance, the graph $\Gamma$ below is chordal, but $\hat\Gamma_2$ is not, since it
has $[v_1,v_3,v_5,v_7]$ as a 4-minimal cycle.

\begin{figure}[ht]
\psfrag{1}{$v_1$}\psfrag{2}{$v_2$}\psfrag{3}{$v_3$}\psfrag{4}{$v_4$}
\psfrag{5}{$v_5$}\psfrag{6}{$v_6$}\psfrag{7}{$v_7$}\psfrag{8}{$v_8$}
\includegraphics[width=2in]{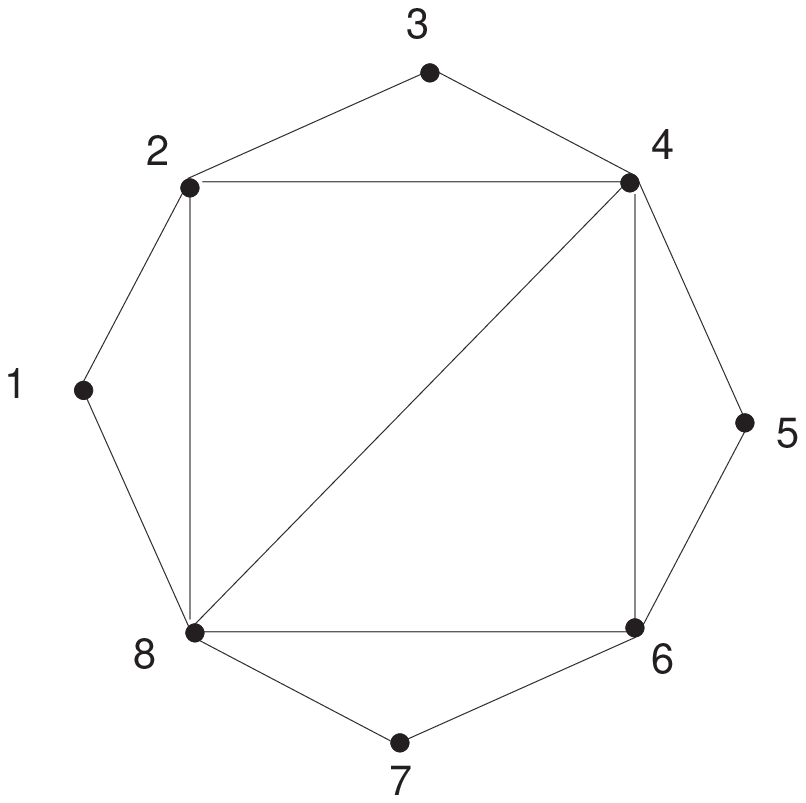}
\label{fi:1}
%	\caption{un test}
\end{figure}

Suppose now that the group $G$ is finitely generated by a set
$A$ with $A=A^{-1}$. The length of an element 
$x\in G$ is defined by
\[
l(x)=\min\{n: x=b_1\cdots b_n, \ b_i\in A  \};
\]
it is equal to the distance between $x$ and $e$ in the Cayley graph
$\Gamma(G,A)$.
If $\Gamma(G,A)$ is a tree, then Lemma~\ref{le:chordal}
and Theorem~\ref{th:main}
yield the following result.

\begin{proposition}\label{pr:tree}
Suppose $G$ is amenable and $\Gamma(G,A)$ is a tree. If $S=\{x\in\Gamma: l(x)\le n\}$, then
any  positive definite function on $S$ can be extended to the whole of~$G$.
\end{proposition}

The proposition applies to
the free product $G=\bbZ_2 \star\bbZ_2$: it is easily seen that,
if $A$ is formed by the two generators, then $\Gamma(G,A)$ is
order isomorphic to $\bbZ$, and is thus a tree. So  
any positive definite function defined on words of length 
smaller or equal to $n$ extends to the whole group.

Unfortunately,
there seem not to be many amenable graphs whose Cayley graph with respect
to some set of generators is a tree. 
Note first the following simple lemma.

\begin{lemma}\label{le:11}
Suppose $G$ is a group, $A\subset G$ is a set of generators, and $\Ga(G,A)$ is a tree. 

(i) For every $x\in G$, 
there is a unique way of writing $x=a_1 \cdots a_n$, with $a_i\in A$, and $a_ia_{i+1}\not=e$; moreover, $l(x)=n$.
(We call $a_1, a_2,...,a_n$ the letters of $x$.)

(ii) Take $x\in G$, with $a_x$ the first letter of $x$. If $y\in G$, and the last letter of $y$ is not $a_x^{-1}$, then $l(yx)=l(x)+l(y)$.
\end{lemma}

We can then obtain the following proposition.

\begin{proposition}\label{pr:tree2}
Suppose that 
$G$ is a discrete amenable group, and $A\subset G$ is a subset of generators, such that 
$\Gamma(G,A)$ is a tree. Then
 either $G=\bbZ$, or $G=\bbZ_2 \star\bbZ_2$.
\end{proposition}

\begin{proof}
Note first that $G$ cannot be finite, since then we may take an element 
$a\in A$ with finite order $p$, and construct the cycle $[e, a, a^2,\ldots, a^{p-1}]$ in $\Ga(G, A)$, which has no chords.

One of the alternate definitions of an amenable group is 
the F\o lner condition, which in the case of discrete groups can be stated as follows: given any finite set 
$F\subset G$ and any $\epsilon>0$, there exists a finite subset $K\subset G$, such that  \[\frac{{\rm card}
 (K \bigtriangleup FK)}{{\rm card} K}<\epsilon\]
($K \bigtriangleup FK$ is the symmetric difference).
Using a translation, if necessary, we may assume $e\in K$. Denote also $S_n=\{x\in G: l(x)=n\}$.

Suppose that $x\in G$; Lemma~\ref{le:11} implies that there is at most one element  $a\in A$ with 
the property that $l(ax)\not= l(x)+1$ (otherwise there would exist a cycle in $\Ga(G,A)$). Therefore, if $x\in S_n$, there is at 
most one $a\in A$ such that $ax\not\in S_{n+1}$. Moreover, if $x,y\in S_n$, $x\not=y$, $a,b\in A$ with $ax, 
by\in S_{n+1}$, then $ax\not=by$ (otherwise we obtain again a cycle in $\Ga(G,A)$.

It follows then that, if $A$ has at least 3 elements, then, for any finite set $E\subset S_n$, $AE\cap S_{n+1}$ 
has at least twice more elements than $E$. Therefore
\begin{equation}\label{eq:55}
{\rm card} K=\sum_n {\rm card}(K\cap S_n)\le 2 \sum_n {\rm card}(AK\cap S_{n+1})\le 2{\rm card}(AK).
\end{equation}
Thus ${\rm card}(K\bigtriangleup AK)\ge {\rm card} K$, and the F\o lner condition cannot be satisfied.

Therefore $A$ has at most two elements. If it has only one element, then, being infinite, it is $\bbZ$. 

Suppose it has two elements. If $a^2\not=e$ and $x\in G$, then applying again Lemma~\ref{le:11}, we have that
$l(a'x)\not=l(x)+2$  
for at most one element $a'$ 
in the set $A'=\{a^2, ab, ba\}$, and for $x,y\in S_n$, $x\not=y$,  $a', b'\in A'$  with 
$a'x, b'y\in S_{n+2}$, we have $a'x\not= b'y$. Therefore, for any  finite set $E\subset S_n$, $AE\cap S_{n+2}$ has at 
least twice more elements than $E$, and we obtain~\eqref{eq:55} with $S_{n+1}$ replaced by $S_{n+2}$. Thus again 
${\rm card}(K\bigtriangleup AK)\ge {\rm card} K$, and the F\o lner condition cannot be satisfied.

Since a similar argument applies in case $b^2\not=e$, the only remaining possibility is $a^2=b^2=e$. Now if either $ab$ or $ba$ 
would have finite order, this would produce a cycle in $\Ga(G, A)$. Thus they are both of infinite order, and it follows easily 
that $G$ is isomorphic to $\bbZ_2 \star\bbZ_2$.
\end{proof}

\end{document}